\documentclass[graybox]{svmult}

\usepackage{mathptmx}       
\usepackage{amssymb,amsmath,mathtools}
\usepackage{helvet}         
\usepackage{courier}        
\usepackage{makeidx}         
\usepackage{graphicx}        
\usepackage{multicol}        
\usepackage[bottom]{footmisc}
\usepackage{array,booktabs}
\newcolumntype{C}[1]{>{\centering\arraybackslash}m{#1}}

\makeindex             

\newtheorem{thm}{Theorem}

\begin{document}

\title*{Random Noble Means Substitutions}
\author{Michael Baake and Markus Moll}
\institute{Fakult\"{a}t f\"{u}r Mathematik, Universit\"{a}t Bielefeld, 
Postfach 100131, 33501 Bielefeld, Germany\newline \and 
Michael Baake; \email{mbaake@math.uni-bielefeld.de}
\newline \and Markus Moll; \email{mmoll@math.uni-bielefeld.de}}
%
%
\maketitle

\abstract*{The random local mixture of a family of primitive substitution rules with
noble mean inflation multiplier is investigated. This extends the random
Fibonacci example that was introduced by Godr\`{e}che and Luck in 1989. We
discuss the structure of the corresponding dynamical systems, and determine their
entropy, an ergodic invariant measure and diffraction spectrum.}

\abstract{The random local mixture of a family of primitive substitution rules with
noble mean inflation multiplier is investigated. This extends the random
Fibonacci example that was introduced by Godr\`{e}che and Luck in 1989. We
discuss the structure of the corresponding dynamical systems, and determine their
entropy, an ergodic invariant measure and diffraction spectrum.}

\section{Introduction}\label{sec:1}

Despite many open problems (including the famous Pisot substitution conjecture), the
structure of systems with pure point diffraction is rather well understood
\cite{bamo_bamoo,bamo_queff}.  Due to recent progress \cite{bamo_babimo,
bamo_bagagr}, also the situation for various systems with diffraction spectra of
mixed type has improved, in particular from a computational point of view. In
particular, one can explicitly calculate the diffraction measure in closed form for
certain classes of examples.  Still, the understanding of spectra in the presence of
entropy is only at its beginning; compare \cite{bamo_babimo,bamo_baleri} and
references therein. The purpose of this contribution is a further step into
'disordered territory', here via the analysis of mixed substitutions that are
randomly applied at a \emph{local} level. This is in contrast to global mixtures
(which leads to $S$-adic systems), which have no entropy.  Local mixtures were
investigated in \cite{bamo_luck}, where the essential properties of the Fibonacci
case were derived, along with first results on planar systems based on triangle
inflation rules. Here, we extend the random Fibonacci system to the noble means
family, and present the results from the point of view of dynamical systems.  The
entire family is still relatively simple because each individual member of a fixed
noble mean family defines the same (deterministic) hull.  Various generalisations are
possible, but not discussed here.

\section{Construction}\label{sec:2}
Let $\mathcal{A} = \{a,b\}$ be our two letter alphabet. For any fixed integer $m \ge
1$, we define a family $\mathcal{H}_{m}$ of substitution rules by
\begin{equation*}
    \zeta_{m,i}^{} \colon
    \left\{
    \begin{array}{lll}
        a \mapsto & a^{i} b a^{m - i} \\
        b \mapsto & a
    \end{array}
    \right.
\end{equation*}
with $0 \le i \le m$, and refer to each $\zeta_{m,i}^{}$ as a \emph{noble means
substitution (NMS)}. Each member of $\mathcal{H}_{m}$ is a primitive substitution
with Pisot inflation multiplier $\lambda_{m} = \frac{1}{2}(m + \sqrt{m^{2} + 4})$ and
algebraic conjugate $\lambda_{m}^{\prime} = \frac{1}{2}(m - \sqrt{m^{2} + 4})$. Each
substitution possesses a reflection symmetric and aperiodic two-sided discrete (or
symbolic) hull $\mathbb{X}_{m,i}^{}$, where the hull, as usual, is defined as the
orbit closure of a fixed point in the local topology. Moreover, all elements of
$\mathcal{H}_{m}$ are pairwise conjugate to each other which implies that, for fixed
$m$, the hulls $\mathbb{X}_{m,i}^{}$ are equal for $0 \le i \le m$.

We now fix a probability vector $(p_{0}^{}, \ldots ,p_{m}^{})$, that is $p_{i}^{} \ge
0$ and $\sum_{j = 0}^{m} p_{j}^{} = 1$. We define the random substitution rule
\begin{equation*}
    \zeta_{m}^{} \colon 
    \left\{
    \begin{array}{lll}
        a  \mapsto & \left\{
        \begin{array}{cc}
            \zeta_{m,0}^{}(a) \,=\, ba^{m}, &\text{with probability } p_{0}^{}, \\
            \vdots                          &\vdots                             \\
            \zeta_{m,m}^{}(a) \,=\, a^{m}b, &\text{with probability } p_{m}^{},
        \end{array}
        \right. \\
        b \mapsto & a
    \end{array}
    \right.
    \quad \text{where} \quad
    M \,:=\,
    \begin{pmatrix}
        m & 1 \\
        1 & 0
    \end{pmatrix}
\end{equation*}
is its substitution matrix. We refer to $\zeta_{m}^{}$ as the \emph{random noble
means substitution (RNMS)}.  The application of $\zeta_{m}^{}$ occurs \emph{locally},
which means that we decide separately on each letter $a$ which of the $m + 1$
possible realisations we choose.  In particular, for each $k \in \mathbb{N}$,
$\zeta_{m}^{k}(a)$ is a random variable. As there is no direct analogue of a fixed
point in the stochastic situation, we have to slightly adjust the notion of the
two-sided discrete hull in this context. Note that $aa$ is a legal word (see below
for more) for all $m$, and consider
\begin{equation*}
    X_{m}^{} \,:=\, \bigl\{w \in \mathcal{A}^{\mathbb{Z}} \mid w \text{ is an
    accumulation point of } \bigl( S^{j_{n}} \zeta_{m}^{n}(a | a) \bigr)_{n \in
    \mathbb{N}_{0}} \bigr\},
\end{equation*}
where $S$ denotes the shift.  The two-sided discrete hull $\mathbb{X}_{m}^{}$ is
defined as the smallest closed and shift-invariant subset of
$\mathcal{A}^{\mathbb{Z}}$ with $X_{m}^{} \subset \mathbb{X}_{m}^{}$. It is immediate
that $\mathbb{X}_{m}^{}$ is a superset of $\mathbb{X}_{m,i}^{}$. Note that typical
elements of $\mathbb{X}_{m}^{}$ contain the subword $bb$, which is absent in
$\mathbb{X}_{m,i}^{}$. The hull $\mathbb{X}_{m}^{}$ is characterised by the property
that it contains all elements of $\{ a,b \}^{\mathbb{Z}}$ that contain
$\zeta_{m}^{}$-legal subwords only (see below for more).

\section{Topological entropy}\label{sec:3}

In this section, we assume that all probabilities $p_{i}^{}$ are strictly positive.
We call a finite word $w$ \emph{legal} with respect to $\zeta_{m}^{}$ if there is a
power $k \in \mathbb{N}$ such that $w$ is a subword of some realisation of
$\zeta_{m}^{k}(a)$.  Furthermore, let $\mathcal{D}_{\!m,\ell}^{}$ be the set of all
legal words of length $\ell$ with respect to $\zeta_{m}^{}$. We refer to the function
$C \colon \mathbb{N} \to \mathbb{N}$, $\ell \mapsto \lvert \mathcal{D}_{\!m,\ell}^{}
\rvert$ as the \emph{complexity function} of $\zeta_{m}^{}$. It is known that the
discrete hull of each member of $\mathcal{H}_{m}$ has linear complexity, which
implies that the topological entropy vanishes here. In the stochastic setting, the
picture changes; see \cite{bamo_baleri} for background.

Let $m \in \mathbb{N}$ be arbitrary but fixed. The sets $\mathcal{G}_{0}^{} :=
\varnothing$, $\mathcal{G}_{1}^{} := \{b\}$, $\mathcal{G}_{2}^{} := \{a\}$ and
\begin{equation}\label{equ:1}
    \mathcal{G}_{n}^{} \,:=\, \bigcup_{i = 0}^{m} \prod_{j = 0}^{m} \mathcal{G}_{n -
        1 - \delta_{ij}}^{},
\end{equation}
with $\delta_{ij}^{}$ denoting the Kronecker symbol, are called the \emph{generation sets}
of $\zeta_{m}^{}$. The product in $\eqref{equ:1}$ is meant to be the set-theoretic product
with respect to concatenation of words. Moreover, we define $\mathcal{G} := \lim_{n \to
\infty} \mathcal{G}_{n}^{}$ and refer to $\mathcal{G}_{n}^{}$ as the \emph{$n$-th
generation set}. The length $\ell_{n}^{}$ of words in $\mathcal{G}_{n}^{}$ is given by the
sequence $\ell_{0}^{} = 0$, $\ell_{1}^{} = 1$, $\ell_{2} = 1$ and $\ell_{n + 1}^{} =
m\ell_{n}^{} + \ell_{n - 1}^{}$, for $n \ge 2$. The set $\mathcal{G}_{n}^{}$ consists of
all possible exact realisations $\zeta_{m}^{n-1}(b)$. Since not all legal words result
from an exact substitution, which can again be seen from the example $bb$, it is clear
that $\lvert \mathcal{G}_{\!n}^{} \rvert < C(\ell_{n}^{})$ for $n \ge 2$.

In \cite{bamo_luck}, Godr\`{e}che and Luck computed the topological entropy of
$\zeta_{1}^{}$ under the implicit assumption that
\begin{equation*}
    \lim_{n \to \infty} \frac{1}{\ell_{n}^{}}\log(C(\ell_{n}^{})) \,=\, \lim_{n \to
    \infty} \frac{1}{\ell_{n}^{}} \log(\lvert \mathcal{G}_{n}^{} \rvert),
\end{equation*}
which was recently proved by J.\ Nilsson \cite{bamo_nilsson}. This asymptotic
identity is crucial because the exact computation of the complexity function of
$\zeta_{m}^{}$ is still an open problem. It is easy to compute $\lvert
\mathcal{D}_{\!m,\ell}^{} \rvert$ for $\ell \le m + 2$ and it is known
\cite{bamo_moll} that
\begin{equation*}
    \lvert \mathcal{D}_{\!m,\ell}^{} \rvert \,=\, \sum_{i = 0}^{3} \binom{\ell}{i} -
    \frac{1}{6}m(m + 1)(3\ell - 2m - 4)
\end{equation*}
if $m + 3 \le \ell \le 2m + 2$, while an extension to arbitrary word lengths seems
difficult. In \cite{bamo_luck}, the entropy per letter for $m = 1$ is computed to be
\begin{equation*}
    h_{1}^{} \,=\, \lim_{n \to \infty} \frac{\log(\lvert \mathcal{G}_{n}^{}
    \rvert)}{\ell_{n}^{}} \,=\, \sum_{i = 2}^{\infty} \frac{\log(i)}{\lambda_{1}^{i+2}}
    \,\approx\, 0.444399 \,>\, 0,
\end{equation*}
whereas a convenient representation for arbitrary $m$ reads
\begin{equation}\label{equ:2}
    h_{m}^{} \,=\, \lim_{n \to \infty} \frac{\log(\lvert \mathcal{G}_{n}^{}
    \rvert)}{\ell_{n}^{}} \,=\, \frac{\lambda_{m}^{} - 1}{1 - \lambda_{m}^{\prime}} \cdot
    \sum_{i = 2}^{\infty} \frac{\log(m(i - 1) + 1)}{\lambda_{m}^{i}}.
\end{equation}
The result computed by Godr\`{e}che and Luck for $m = 1$ can be recovered by the
observation that $(\lambda_{1}^{} - 1)/(1 - \lambda_{1}^{\prime}) =
1/\lambda_{1}^{2}$ in this case. Some numerical values are given in Table
$\ref{TAB:1}$.
\begin{table}
    \centering
    \begin{tabular}{C{1.25cm}C{1.25cm}C{1.25cm}C{1.25cm}C
        {1.25cm}C{1.25cm}C{1.25cm}C{1.25cm}}
        \toprule
        $m$  & $1$ & $2$ & $3$ & $4$ & $5$ & $6$ & $7$ \\
        \midrule
        $h_{m}^{}$ & $0.444399$ & $0.408549$ & $0.371399$ & $0.338619$ &
        $0.310804$ & $0.287298$ & $0.267301$ \\
        \bottomrule
    \end{tabular}
    \caption{Numerical values of the topological entropy for RNMS with $1 \le m
            \le 7$.\label{TAB:1}}
\end{table}
It is not difficult to prove \cite{bamo_moll} that $\lim_{m \to \infty} h_{m}^{} =
0$, which can be verified by estimating the logarithm in $\eqref{equ:2}$ via the
square root and using the fact that $\lambda_{m}^{}/m$ tends to $1$ as $m \to
\infty$.

\section{Frequencies of subwords}\label{sec:4}

We adopt the method of computing the frequencies of subwords via induced
substitutions on words of length $\ell$ (with $\ell \in \mathbb{N}$), which was
introduced in \cite[Section 5.4.1]{bamo_queff}, and modify it to fit the stochastic
setting. To this end, we again assume that all probabilities $p_{i}^{}$ in the
definition of $\zeta_{m}^{}$ are strictly positive.

If $w = w_{0}^{} w_{1}^{} \cdots w_{\ell - 1}^{}$ is a word of length $\ell$, we
define $w_{[i,j]}^{}$ to be the subword $w_{i}^{} \cdots w_{j}^{}$ of $w$ of length
$j - i + 1$.  For $\ell \ge 2$, we denote $\zeta_{m,(\ell)}^{} \colon
\mathcal{D}_{\!m,\ell}^{} \to \mathcal{D}_{\!m,\ell}^{}$ as the \emph{induced
substitution} defined by
\begin{equation}\label{equ:induced}
    \zeta_{m,(\ell)}^{} \colon w^{(i)} \mapsto
    \left\{
    \begin{array}{cc}
        \left( v^{(i,1)}_{[k,k + \ell - 1]} \right)_{0 \le k \le
        \lvert \zeta_{m}(w_{0}^{(i)}) \rvert - 1}, &\text{with probability } p_{i,1}^{},
        \\
        \vdots & \vdots \\
        \left( v^{(i,n_i)}_{[k,k + \ell - 1]} \right)_{0 \le k \le
        \lvert \zeta_{m}(w_{0}^{(i)}) \rvert - 1}, &\text{with probability }
        p_{i,n_i}^{}, \\
    \end{array}
    \right.
\end{equation}
where $w^{(i)} \in \mathcal{D}_{\!m,\ell}^{}$ and $v^{(i,j)}$ is a realisation of
$w^{(i)}$ under $\zeta_{m}^{}$ with probability $p_{i,j}^{}$. This way, we ensure
that we are neither under- nor overcounting subwords of a given length. Similar to
the case of $\zeta_{m}^{}$, the result is a random variable.

The action of $\zeta_{m,(\ell)}^{}$ on words in $\mathcal{D}_{\!m,\ell}^{}$ is
illustrated in the following table for $m = 1$ and $\ell = 2$:
\begin{center}
    \begin{tabular}{ccccccccc}
        \toprule
        $w \in \mathcal{D}_{\!1,2}^{}$ & $\zeta_{1}^{}(w)$ & $v^{(i,j)}$ &
        $\mathbb{P}$ & \phantom{\qquad}& $w \in \mathcal{D}_{\!1,2}^{}$ & $\zeta_{1}^{}(w)$ &
        $v^{(i,j)}$ & $\mathbb{P}$ \\
        \midrule
        $aa$ & $abab$ & $(ab)(ba)$ & $p_{1}^{2}$ && $ab$ & $aba$ & $(ab)(ba)$ & $p_{1}^{}$ \\
             & $abba$ & $(ab)(bb)$ & $p_{0}^{}p_{1}^{}$ &&      & $baa$ & $(ba)(aa)$ & $p_{0}^{}$ \\
             & $baab$ & $(ba)(aa)$ & $p_{0}^{}p_{1}^{}$ && $ba$ & $aab$  & $(aa)$ & $p_{1}^{}$  \\
             & $baba$ & $(ba)(ab)$ & $p_{0}^{2}$  && & $aba$  & $(ab)$     & $p_{0}^{}$         \\
             &&&&&$bb$ & $aa$   & $(aa)$     & $1$                \\
        \bottomrule
    \end{tabular}
\end{center}

\noindent Applying the lexicographic order to the words in
$\mathcal{D}_{\!m,\ell}^{}$ leads to the corresponding substitution matrix
$M_{m,\ell}^{} := M(\zeta_{m,(\ell)}^{})$. For any fixed $m \in \mathbb{N}$ and $\ell
= 2$, we get
\begin{equation*}
    M_{m,2}^{} \,=\,
    \begin{pmatrix*}
        (m-1)+p_{0}^{}p_{m}^{} & (m-1)+p_{0}^{} & 1-p_{0}^{} & 1 \\
        1-p_{0}^{}p_{m}^{}     & 1-p_{0}^{}     & p_{0}^{}   & 0 \\
        1-p_{0}^{}p_{m}^{}     & 1              & 0          & 0 \\
        p_{0}^{}p_{m}^{}       & 0              & 0          & 0
    \end{pmatrix*}.
\end{equation*}
This matrix has the spectrum $\sigma \bigl(M_{m,2}^{}\bigr) =
\{\lambda_{m}^{},\lambda_{m}^{\prime},-p_{0}^{},p_{0}^{}p_{m}^{} \}$.  Furthermore,
it is interesting to observe that the spectrum of the matrices $M_{m,\ell}$ for $\ell
\ge 3$ is the same as that of $M_{m,2}^{}$, except for the addition of zeros. Note
that $\zeta_{m,(1)}^{}$ agrees with $\zeta_{m}^{}$ which implies that $M_{m,1}^{} =
M$.

The substitution matrix $M_{m,\ell}^{}$ is primitive for all $m$ and $\ell$, which
allows an application of Perron-Frobenius theory; see \cite{bamo_seneta} for general
background.  This implies that there is a strictly positive right eigenvector
$\phi^{(\ell)}$ to the eigenvalue $\lambda_{m}^{}$. Note that $\lambda_{m}^{}$ does
not depend on any of the $p_{i}^{}$ at all, whereas this is not the case for
$\phi^{(\ell)}$.

We define a measure on the discrete hull $\mathbb{X}_{m}^{}$ as follows. For any word
$v \in \mathcal{D}_{\!m,\ell}^{}$ and $k \in \mathbb{N}$, let $Z_{k}^{}(v) :=
\bigl\{w \in \mathbb{X}_{m}^{} \mid w_{[k,k + \ell - 1]}^{} = v \bigr\}$ be the
cylinder set of $v$ that starts at position $k$.  Then, the family
$\{Z_{k}^{}(v)\}_{k \in \mathbb{N}}$ generates the product topology and we define the
measure $\mu \colon \mathbb{X}_{m}^{} \to \mathbb{R}_{\ge 0}^{}$ on the cylinder sets
as $\mu\bigl(Z_{k}^{}(v)\bigr) = \phi^{(\ell)}(v)$, where $\phi^{(\ell)}(v)$ is the
entry of $\phi^{(\ell)}$ with respect to $v$. This is a proper (and consistent)
definition of a measure on $\mathbb{X}_{m}^{}$, which can also be found in
\cite[Section 5.4.2]{bamo_queff}. By construction, the measure is shift-invariant.

The following theorem \cite{bamo_moll} shows that, similar to the deterministic
setting \cite{bamo_queff}, it is possible to interpret the entries of $\phi^{(\ell)}$
as the frequencies of legal subwords with respect to $\zeta_{m}^{}$ as follows:

\begin{thm}\label{THM:ergodic}
    Let\, $\mathbb{X}_{m}^{} \subset \mathcal{A}^{\mathbb{Z}}$ be the hull of the
    random noble means substitution for $m \in \mathbb{N}$ and $\mu$ the
    shift-invariant probability measure on $\mathbb{X}_{m}^{}$ as defined above. For
    any $f \in L^{1}(\mathbb{X}_{m}^{},\mu)$ and for an arbitrary but fixed $s \in
    \mathbb{Z}$, the identity
    \begin{equation*}
        \lim_{N \to \infty} \frac{1}{N} \sum_{i = s}^{N + s - 1} f(S^{i}x) \,=\,
        \int_{\mathbb{X}_{m}^{}}f \,\mathrm{d} \mu
    \end{equation*}
    holds for $\mu$-almost every $x \in \mathbb{X}_{m}^{}$.
\end{thm}

The proof can be accomplished by inspecting the family of random variables $R =
\{f(S^{i}w)\}_{i \in \mathbb{N}}$, where $f$ is a \emph{patch recognition function}
that evaluates to $1$ if $(S^{i}w)_{[s,s + \ell - 1]} = v$ for an arbitrary but fixed
word $v \in \mathcal{D}_{\!m,\ell}^{}$ and to $0$ otherwise. By observing that, given
any $i \in \mathbb{Z}$, the sets
\begin{equation*}
    \mathcal{I}_{i} \,:=\, \bigl\{\,(S^{i + k(\ell + m)}w)_{[s,s + \ell -
    1]}^{} \mid k \in \mathbb{N}\,\bigr\}
\end{equation*}
consist of pairwise independent words, we can split up the summation over $R$
appropriately and apply Etemadi's version of the strong law of large numbers
\cite[Theorem 1]{bamo_etemadi} to each sum over $\mathcal{I}_{i}^{}$ separately.
This, in conjunction with an application of the Stone-Weierstrass theorem, implies
the assertion.

\section{Diffraction measure}\label{sec:5}

The symbolic situation is turned into a geometric one as follows. In view of the left
PF eigenvector of $M$, $a$ and $b$ are turned into intervals of lengths
$\lambda_{m}^{}$ and $1$, respectively.  The left end points are the coordinates we
use. The corresponding \emph{continuous} hull $\mathbb{Y}_{\!m}^{}$ is the orbit
closure of all accumulation points of the geometric inflation rule under
$\mathbb{R}$. Let $\Lambda \subset \mathbb{Z}[\lambda_{m}^{}]$ be a coordinatisation
of an element of $\mathbb{X}_{m}^{}$ in $\mathbb{R}$.  Then, $\Lambda$ is a discrete
point set that fits into the same cut and project scheme as all elements of the
family $\mathcal{H}_{m}^{}$.  With respect to $\Lambda$, the smallest interval that
covers $\Lambda^{\prime} = \{x^{\,\prime} \mid x \in \Lambda \}$ in internal space is
given by $[\lambda_{m}^{\prime} - 1,1 - \lambda_{m}^{\prime}]$.  Then, $\Lambda$ is
relatively dense with covering radius $\lambda_{m}^{}$, and a subset of a model set,
which implies that $\Lambda$ is a Meyer set by \cite[Theorem 9.1]{bamo_moody}.  Let
$\Lambda_{R}^{} = \Lambda \cap \overline{B_{R}^{}}$ and consider the autocorrelation
\begin{equation*}
    \gamma \,:=\, \lim_{R \to \infty} \frac{\delta_{R}^{} \ast
    \widetilde{\delta_{R}^{}}}{\mathrm{vol}(B_{R}^{})} \quad \text{with} \quad
    \delta_{R}^{} \,=\, \sum_{x \in \Lambda_{R}^{}} \delta_{x}^{}.
\end{equation*}
The limit almost surely exists due to the ergodicity of our system.  By construction,
$\gamma$ is a positive definite measure which implies that its Fourier transform
exists and is a positive measure. Regarding the Lebesgue decomposition
$\widehat{\gamma} = (\widehat{\gamma})_{\mathrm{pp}}^{} +
(\widehat{\gamma})_{\mathrm{ac}}^{} + (\widehat{\gamma})_{\mathrm{sc}}^{}$, it is
possible to compute the pure point part to be
\begin{equation*}
    (\widehat{\gamma})_{\mathrm{pp}}^{} \,=\, \sum_{k \in L^{\circledast}} \lvert
    \widehat{\eta}_{a}(-k^{\prime}) + \widehat{\eta}_{b}(-k^{\prime}) \rvert^{2}
    \delta_{k}^{},
\end{equation*}
where $L^{\circledast} = \mathbb{Z}[\lambda_{m}^{}]/\sqrt{m^{2} + 4}$  is the Fourier
module. In the case of $m=1$, the invariant measures $\widehat{\eta}_{a}^{}$,
$\widehat{\eta}_{b}^{}$ can be approximated via the recursion relation
\begin{equation}\label{equ:4}
    \begin{pmatrix}
        \widehat{\eta}_{a}^{}(y) \\
        \widehat{\eta}_{b}^{}(y)
    \end{pmatrix}
    \,=\,
    \lvert \xi \rvert^{n} \cdot \prod\limits_{\ell=1}^{n} 
    \left[
    p_{0}^{}
    \begin{pmatrix*}[l]
        \mathrm{e}^{-2\pi\mathrm{i}y\xi^{\ell-1}} & 1 \\
        1                                         & 0
    \end{pmatrix*}
    + p_{1}^{}
    \begin{pmatrix*}[l]
        1                                       & 1 \\
        \mathrm{e}^{-2\pi\mathrm{i}y\xi^{\ell}} & 0
    \end{pmatrix*} 
    \right]
    \cdot
    \begin{pmatrix}
        \widehat{\eta}_{a}^{}(y\xi^{n}) \\
        \widehat{\eta}_{b}^{}(y\xi^{n})
    \end{pmatrix},
\end{equation}
with $n \in \mathbb{N}$ and $\xi := \lambda_{1}^{\prime}$. As $\xi^{n} \to 0$ for $n
\to \infty$, an appropriate choice of the eigenvector
$\bigl(\widehat{\eta}_{a}^{}(0),\widehat{\eta}_{b}^{}(0)\bigr)^{T}$ for the equation
\begin{equation*}
    \begin{pmatrix}
        1 & 1 \\
        1 & 0
    \end{pmatrix}
    \cdot
    \begin{pmatrix}
        \widehat{\eta}_{a}^{}(0) \\
        \widehat{\eta}_{b}^{}(0)
    \end{pmatrix}
    =
    \lambda_{1}^{} \cdot
    \begin{pmatrix}
        \widehat{\eta}_{a}^{}(0) \\
        \widehat{\eta}_{b}^{}(0)
    \end{pmatrix},
\end{equation*}
which results from $\eqref{equ:4}$ for $k = 0$ and $n = 1$, fixes the base of the
recursion and provides the desired approximation. Since $\widehat{\eta}_{a}(0) +
\widehat{\eta}_{b}(0)$ must be the density of $\Lambda$, which always is
$\lambda_{1}^{}/\sqrt{5}$, one finds $\widehat{\eta}_{a}(0) = 1/\sqrt{5}$ and
$\widehat{\eta}_{b}(0) = (\lambda_{1}^{}-1)/\sqrt{5}$. Let $\mu_{s}^{}$ be the
measure on $\mathbb{Y}_{\!m}^{}$ induced by $\mu$ via suspension; see \cite[Chapter
11]{bamo_sinai} for general background. One consequence, due to a theorem of
Strungaru \cite{bamo_strungaru} and an application of the methods of
\cite{bamo_bale}, is that our random dynamical system $D :=
(\mathbb{Y}_{\!m}^{},\mathbb{R},\mu_{s}^{})$ is ergodic, but not weakly mixing. In
particular, it has strong long-range order.

Due to the stochastic setting with positive entropy, one expects a non-trivial
absolutely continuous part. For $m=1$, the ergodicity of $D$ almost surely yields a
diffraction measure which can be represented as $\widehat{\gamma} \,=\,
(\widehat{\gamma})_{\mathrm{pp}}^{} + \alpha(k) \cdot \lambda$, where
\begin{equation*}\label{equ:3}
    \widehat{\gamma}(\{k\}) \,=\, \lim_{n \to \infty} \frac{1}{L_{n}^{2}} \cdot
    \bigl\lvert \mathbb{E} \bigl(g_{n}^{}(k) \bigr) \bigr\rvert^{2}.
\end{equation*}
Here, $\mathbb{E}$ refers to averaging with respect to $\mu_{s}^{}$ and
\begin{equation*}
    \alpha(k) \,:=\, \lim_{n \to \infty} \frac{1}{L_{n}^{}} \cdot \left(
    \mathbb{E}\bigl( \bigl\lvert g_{n}^{}(k) \bigr\rvert^{2} \bigr) - \bigl\lvert
    \mathbb{E}\bigl( g_{n}^{}(k) \bigr) \bigr\rvert^{2} \right),
\end{equation*}
with the random exponential sum $g_{n}^{}(k) := \sum_{j=1}^{F_{n+1}}
\mathrm{e}^{-2\pi\mathrm{i} k x_{j}}$ and $L_{n}^{} := \lambda_{1}^{}F_{n}^{} +
F_{n-1}^{}$, where $F_{n}^{}$ is the $n$-th Fibonacci number. Now let
\begin{equation*}
    A_{n}^{}(k) \,:=\, \mathbb{E} \bigl(g_{n}^{}(k)\bigr) \quad \text{and} \quad
    B_{n}^{}(k) \,:=\, \mathbb{E}\bigl( \bigl\lvert g_{n}^{}(k) \bigr\rvert^{2} \bigr)
    - \bigl\lvert \mathbb{E}\bigl( g_{n}^{}(k) \bigr) \bigr\rvert^{2}.
\end{equation*}
Godr\`{e}che and Luck \cite{bamo_luck} derived a recursion relation for the sequence
$A_{n}^{}(k)$,
\begin{equation*}
    A_{n}^{}(k) \,=\, \bigl(p_{1}^{} + p_{0}^{} \mathrm{e}^{-2\pi\mathrm{i}kL_{n-2}}
    \bigr) \cdot A_{n-1}^{}(k) + \bigl(p_{0}^{} + p_{1}^{}
    \mathrm{e}^{-2\pi\mathrm{i}kL_{n-1}} \bigr) \cdot A_{n-2}^{}(k),
\end{equation*}
where $A_{0}^{}(k) = \mathrm{e}^{-2\pi\mathrm{i}k}$ and $A_{1}^{}(k) =
\mathrm{e}^{-2\pi\mathrm{i}k\lambda_{1}}$. Analogously, one derives a recursion
relation for the sequence $B_{n}^{}(k)$,
\begin{equation*}
    B_{n}^{}(k) \,=\, B_{n-1}^{}(k) + B_{n-2}^{}(k) + 2p_{0}^{}p_{1}^{} \cdot
    \Delta_{n}^{}(k),
\end{equation*}
with
\begin{align*}
    \Delta_{n}^{}(k) &\,=\, \bigl(1-\cos(2\pi k L_{n-1}^{}) \bigr) \cdot
    \lvert A_{n-2}^{}(k) \rvert^{2} + \bigl(1-\cos(2\pi k L_{n-2}^{}) \bigr) \cdot
    \lvert A_{n-1}^{}(k) \rvert^{2} \\
    &\qquad - \mathrm{Re} \bigl[ (1 - \mathrm{e}^{2\pi\mathrm{i}kL_{n-1}^{}}) \cdot
    (1 - \mathrm{e}^{-2\pi\mathrm{i}kL_{n-2}^{}}) \cdot A_{n-1}^{}(k) \cdot
    \overline{A_{n-2}^{}(k)} \,\bigr]
\end{align*}
and $B_{0}^{}(k) = B_{1}^{}(k) = 0$. In \cite{bamo_luck}, almost surely by way of a
misprint, the authors applied complex conjugation on $A_{n-1}^{}(k)$ instead of
$A_{n-2}^{}(k)$, which makes a huge difference, as the sequence $B_{n}^{}(k)$ does
not converge in that case. The recursion for $B_{n}^{}(k)$ can be solved and the
explicit representation reads
\begin{equation*}
    B_{n}^{}(k) = 2p_{0}^{}p_{1}^{} \cdot \sum_{i = 2}^{n} F_{n+1-i}^{} \,
    \Delta_{i}^{}(k).
\end{equation*}
A detailed discussion of the continuous part of $\widehat{\gamma}$ can be found in
\cite{bamo_moll}.

\begin{figure}
    \centering
    \includegraphics[scale=0.6]{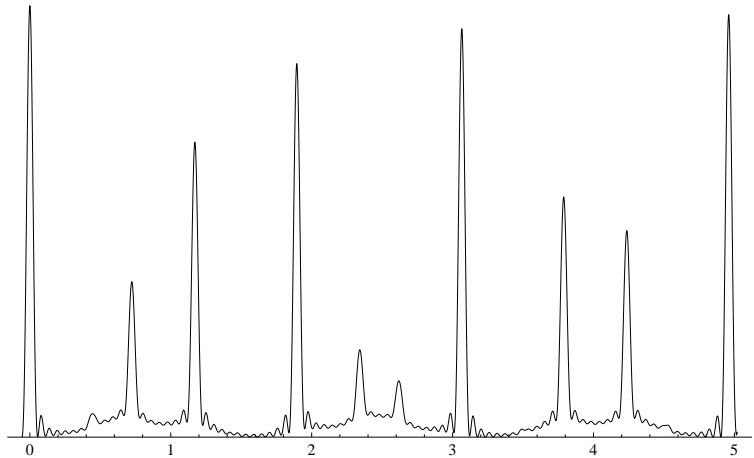}
    \caption{Approximative calculation of the diffraction measure $\widehat{\gamma}$ for
    $m = 1$ and $p_{0}^{} = p_{1}^{} = \frac{1}{2}$, based on $A_{n}^{}(k)$ and
    $B_{n}^{}(k)$ with $n = 6$.\label{fig:1}}
\end{figure}

The illustration of an approximation of the diffraction measure $\widehat{\gamma}$ in
case of $m=1$ and $p_{0}^{} = p_{1}^{} = \frac{1}{2}$, based on the sequences
$A_{n}^{}(k)$ and $B_{n}^{}(k)$, is shown in Figure $\ref{fig:1}$, which agrees with the
average over many realisations for the same length.

\begin{acknowledgement}
    We thank Uwe Grimm and Johan Nilsson for discussions. This work was
    supported by the German Research Council (DFG), within the CRC 701.
\end{acknowledgement}

\end{document}